\documentclass{article}
\usepackage[disable]{todonotes}
\usepackage{amsmath}
\usepackage{graphicx}
\usepackage{caption}
\usepackage{subcaption}
\usepackage{float}
\usepackage{enumitem}
\usepackage{amssymb}
\usepackage[parfill]{parskip}
\usepackage[margin=1.4in]{geometry}
\begin{document}
\author{Weronika Wrzos-Kaminska\footnote{Clare College, Cambridge, Trinity Lane, Cambridge CB2 1TL, UK} }
\title{A simpler winning strategy for Sim}
\maketitle
\begin{abstract}
We give a simple human-playable winning strategy for the second player in the game
of Sim.
\end{abstract}
\section{Introduction}
In the game of Sim, two players, P1 and P2, compete on a complete graph of six vertices ($K_6$). Each of the players has a colour, say red for P1 and blue for P2. The players alternate turns, starting with P1, claiming one of the previously uncoloured edges of the $K_6$ and colouring it in their colour. 
The first player to be forced to form a triangle of their own colour loses. The game of Sim was first introduced by Simmons in 1969 \cite{Simmons}, and has since then attracted a lot of attention. A technical report by Slany \cite{Slany} provides detailed information about Sim and other Ramsey games, as well as further references to literature on Sim. For more information on this and other mis\`ere-type games, see \cite{Leader}.

It is well known that any two edge-colouring of $K_6$ must contain a monochromatic triangle (as the Ramsey number R(3,3) equals 6). Thus, it is impossible for the game of Sim to end in a tie. Consequently, one of the two players must possesses a winning strategy, ie. a strategy that enables the player to win no matter what moves their opponent makes. Computer searches have shown that it is the second player who possesses a winning strategy  for the game of Sim \cite{Simmons}. Mead, Rosa and Huang \cite{Mead} have provided an explicit winning strategy, remarking however that `a simpler (in terms of the rules to be followed) winning strategy is still desirable'. In this paper, we present a different, rather simpler, winning strategy for the game of Sim. 

\section{Some Definitions}
Before describing the strategy, we need to introduce some terminology. 
Firstly, by a position of a game after $k$ moves, we mean a subgraph of $K_6$ with $k$ coloured edges, each either red or blue. 
Given a position, we can partition the edges of $K_6$ as $E(K_6) = R \cup B \cup N,$
where $R$ is the set of red edges (edges that have already been claimed by P1), $B$ is the set of blue edges (edges that have already been claimed by P2) and $N$ is the set of all uncoloured edges. 

By a \textit{P1-allowed move} we mean an uncoloured edge such that colouring it red would not complete a monochromatic triangle. That is, any move that does not lead to an immediate loss for P1 is a P1-allowed move. A \textit{P2-allowed move} is defined analogously. 

Moreover, we will say that a set X of uncoloured edges is a \textit{P1-allowed set} if $X \cup R$ does not contain any triangles. 
Again, we define a \textit{P2-allowed set} in the analogous way. 

Given a position P, We define a \textit{mini-board} of P to be a complete subgraph of the $K_6$ such that: 
\begin{enumerate}[label=(\roman*), nosep]
\item M contains all the coloured edges and at least one P2-allowed move
\item M is minimal with respect to the above condition, that is, no proper subgraph $K \subset M$ satisfies (i).
\end{enumerate}

Clearly, a mini-board exists whenever the position has at least one P2-allowed move. Note that the mini-board, whenever it exists, is unique up to isomorphism. 
\begin{figure}[H]
\center
\label{fig2}
\includegraphics[width =  0.8\textwidth]{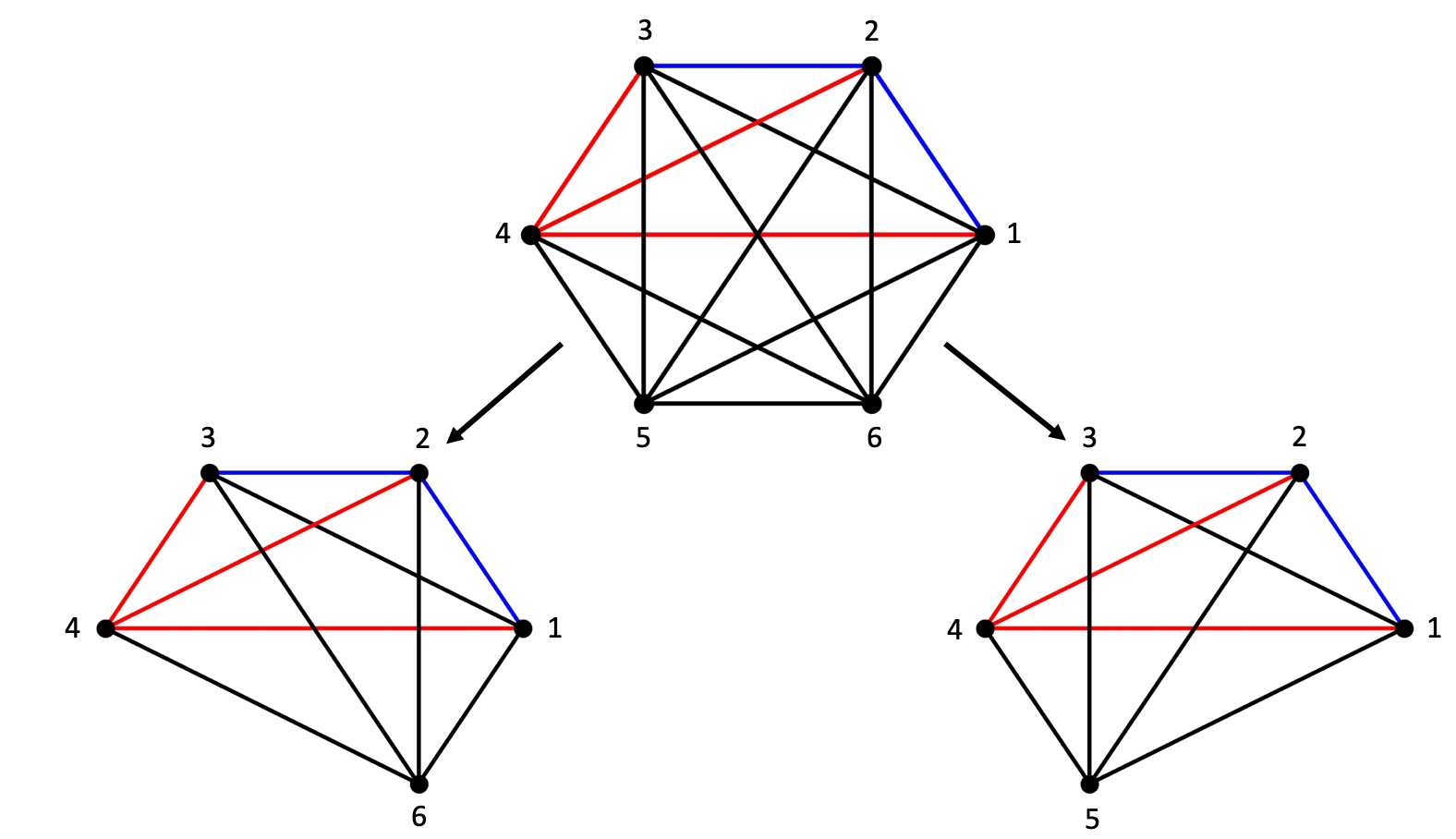}
\caption{Example of a position together with the two (isomorphic) mini-boards. The black edges are representing uncoloured edges. Note that the subgraph spanned by the vertices  $\{1,2,3,4\}$ is \textit{not} a mini-board, as the edge $13$ is not a P2-allowed move.}
\end{figure}
Finally, given a mini-board M, we will say that a set of edges $X \subset E(K_6)$ is a \textit{P1-allowed set on M} if $X \subset E(M)$ and it is a P1-allowed set. By a \textit{maximum P1-allowed set on M} we will mean a P1-allowed set on M of maximum size among all the P1-allowed sets on M.  
As always, we define \textit{P2-allowed set on M} and \textit{maximum P2-allowed set on M} in a similar way. 

\begin{figure}[H]
\center
\label{fig1}
\includegraphics[width =  0.8\textwidth]{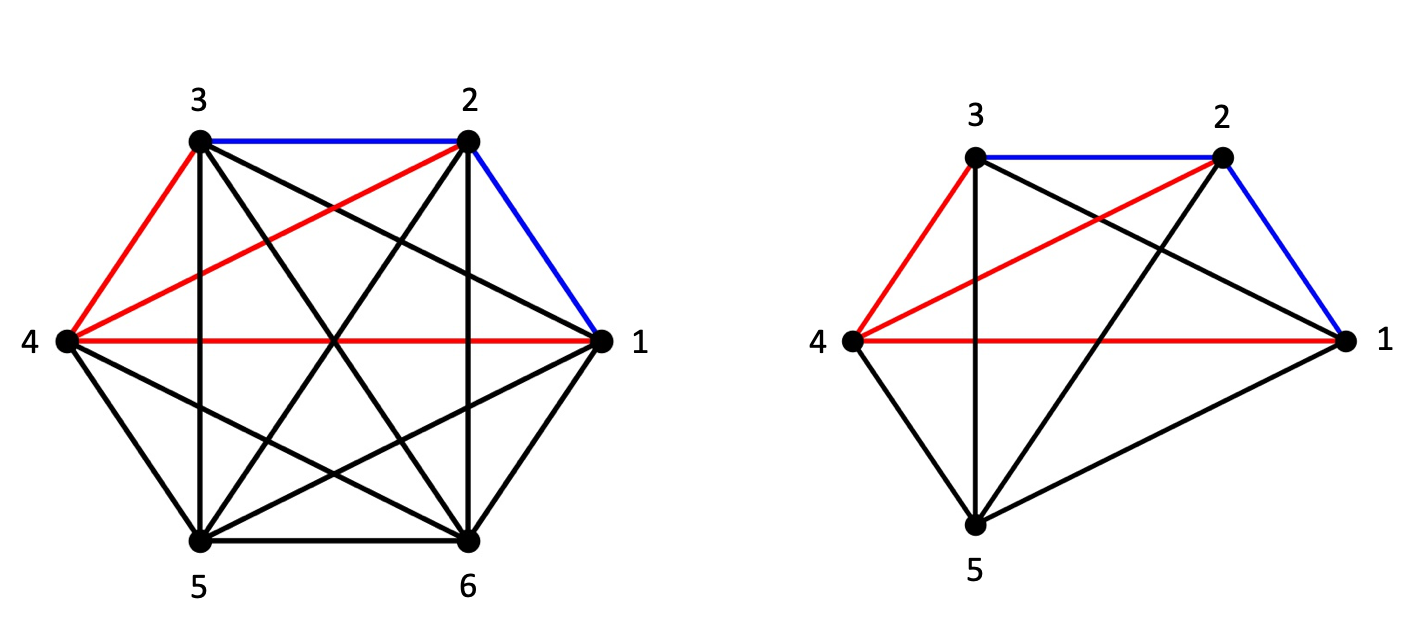}
\caption{For the above position (to the left) and the mini-board M (to the right), we have that the set $\{35, 45 \}$ is a P2-allowed set, but not a P1-allowed set. One can check that there is precisely one maximum P2-allowed set on M, namely $\{15, 35, 45\}$}.
\todo{Fix this caption}
\end{figure}
\section{The strategy}
We are now ready to define the winning strategy for P2: \\
Whenever it is P2's turn, they should fix an arbitrary mini-board M of the current position. They should then pick a move according to the following rules: 
\begin{enumerate}
\item Consider only P2-allowed moves on M. 
\item Pick the move(s) belonging the the greatest number of maximum P2-allowed sets on M. 
\item Pick the move(s)  belonging the the greatest number of maximum P1-allowed sets on M. 
\end{enumerate}
The rules should be interpreted in a hierarchical order as follows:  P2 should start by applying Rule 1. If this does not determine their move uniquely, then Rule 2 should be applied as a `tie-breaker' to distinguish between the moves satisfying Rule 1. If this still does not determine their move, then Rule 3 should be applied as a further tie-breaker. If the move is still not determined uniquely, then P2 may pick arbitrarily among the moves which are left after applying all three rules. 

This strategy does not necessarily determine the move uniquely up to isomorphism. However, if there is more than one move left after applying all the rules, then any of the remaining moves will lead to a win for P2. 

An exhaustive computer search shows that this strategy is indeed winning for P2.

\section{Acknowledgements}
Many thanks go to Professor Imre Leader, my research supervisor, for his valuable support. 
I was provided financial support by the SRIM bursary grant from University of Cambridge, for which I am also grateful. 

\bibliographystyle{acm} 
\bibliography{bibliography.bib}
\end{document}